\documentclass[aos,preprint]{imsart}

\RequirePackage[colorlinks,citecolor=blue,urlcolor=blue]{hyperref}

\usepackage[authoryear]{natbib}
\usepackage{amsmath,amssymb,amsthm,enumerate,epsfig,graphicx}
\usepackage{ifthen,latexsym,syntonly}
\usepackage{natbib}
\usepackage{color}
\usepackage{graphics}
\usepackage{graphicx}
\usepackage{mathrsfs}
\allowdisplaybreaks

\startlocaldefs
\newtheorem{theorem}{Theorem}
\numberwithin{theorem}{subsection}

\numberwithin{corollary}{subsection}

\newtheorem{definition}{Definition}

\newtheorem{lemma}{Lemma}
\numberwithin{lemma}{subsection}
\newtheorem{proposition}{Proposition}
\numberwithin{proposition}{subsection}

\renewcommand{\cite}{\citet*}

\newcommand{\bx}{\mbox{\bf x}}

\newcommand{\bI}{\mbox{\bf I}}

\newcommand{\bX}{\mbox{\bf X}}

\newcommand{\bSigma}{\mbox{\boldmath $\Sigma$}}

\newcommand{\bbP}{\mathbb{P}}
\newcommand{\bPsi}{\mbox{\boldmath $\Psi$}}
\newcommand{\bbeta}{\mbox{\boldmath $\beta$}}
\newcommand{\bbR}{\mathbb{R}}

\newcommand{\argmin}{\mathrm{argmin}}

\endlocaldefs

\begin{document}

\begin{frontmatter}

\title{Nonparametric Screening under Conditional Strictly Convex Loss for Ultrahigh Dimensional Sparse Data}
\runtitle{Goodness of Fit Nonparametric Screening}


\author{\fnms{Xu} \snm{Han}\corref{Xu Han}\ead[label=e1]{hanxu3@temple.edu}}
\address{Department of Statistical Science\\Fox Business School\\Temple University\\Philadelphia, PA 19122\\USA\\ \printead{e1}}
\affiliation{Temple University}

\runauthor{Xu Han}

\begin{abstract}
Sure screening technique has been considered as a powerful tool to handle the ultrahigh dimensional variable selection problems, where the dimensionality $p$ and the sample size $n$ can satisfy the NP dimensionality $\log p=O(n^a)$ for some $a>0$ (Fan \& Lv 2008). The current paper aims to simultaneously tackle the ``universality" and ``effectiveness" of sure screening procedures. For the ``universality", we develop a general and unified framework for nonparametric screening methods from a loss function perspective. Consider a loss function to measure the divergence of the response variable and the underlying nonparametric function of covariates. We newly propose a class of loss functions called conditional strictly convex loss, which contains, but is not limited to, negative log likelihood loss from one-parameter exponential families, exponential loss for binary classification and quantile regression loss. The sure screening property and model selection size control will be established within this class of loss functions. For the ``effectiveness", we focus on a goodness of fit nonparametric screening (Goffins) method under conditional strictly convex loss. Interestingly, we can achieve a better convergence probability of containing the true model compared with related literature. The superior performance of our proposed method has been further demonstrated by extensive simulation studies and some real scientific data example. 
\end{abstract}



\end{frontmatter}


\section{Introduction}
Ultrahigh-dimensional variable selection has become an important problem in modern statistical research due to the big data collection in a variety of scientific areas, such as genomics, bioinformatics, functional magnetic resonance imaging, high frequency finance, etc. In all these problems, statisticians want to select the important covariates associated with the response variable from $p$ covariates. However, the dimensionality $p$ can grow much faster than the sample size $n$. More specifically, $\log p=O(n^a)$ for some $a>0$, which is denoted as nonpolynomial order (NP) (Fan \& Lv 2008). As Fan, Samworth \& Wu (2009) has pointed out: existing variable selection methods based on penalized pseudo likelihood estimation (e.g. Tibshirani 1996, Fan \& Li 2001, Zou \& Hastie 2005, Zou 2006, Cand\`es \& Tao 2007, Zou \& Li 2008, Zhang 2010) can suffer from the simultaneous challenges to computational expediency, statistical accuracy and algorithmic stability in ultrahigh dimensional problems. 

To handle the challenges in the ultrahigh-dimensional problems, Fan \& Lv (2008) introduced a new statistical framework, sure independence screening. Their original method focused on the Gaussian linear regression models, and the important predictors were selected via the marginal correlation ranking.  Formally, let $M_{\star}$ be the set of true important variables, and $\widehat{M}_n$ be the selected variables based on some procedure, then 
\begin{equation}\label{j4}
P(M_{\star}\subset\widehat{M}_n)\geq 1-\epsilon_n
\end{equation}
where $\epsilon_n>0$ and $\epsilon_n\rightarrow0$ as $n\rightarrow\infty$. This is called the ``sure screening property".  Furthermore, the model selection size can be controlled at a polynomial rate of sample size with probability approaching 1.  Because of its powerful performance and computational convenience in ultrahigh dimensional problems, the sure screening framework has received increasing attention in the past few years. Existing literature in this framework have mainly focused on the ``universality" of the screening procedures, that is, developing procedures for various scenarios which possess the ``sure screening property", e.g. generalized linear model by Fan \& Song (2010), nonparametric additive model by Fan, Feng \& Song (2011), rank based model-free feature screening by Zhu, Li, Li \& Zhu (2011), Cox model by Zhao \& Li (2012), robust rank correlation screening by Li, Peng, Zhang \& Zhu (2012), varying coefficient model by Fan, Ma \& Dai (2013),  empirical likelihood based screening by Chang, Tang \& Wu (2013), quantile-adaptive screening by He, Wang \& Hong (2013), censored rank independence screening by Song, Lu, Ma \& Jeng (2014), fused Kolmogorov filter by Mai \& Zou (2015). On the other hand, formal pursuit of ``effectiveness" of sure screening procedures have been largely ignored. Intuitively, for the convergence probability $1-\epsilon_n$ in (\ref{j4}), if $\epsilon_n$ converges to 0 slower, the corresponding screening procedure will have larger possibility of not selecting the true important variables. More specifically, existing literature commonly show that 
\begin{equation}\label{j5}
P(M_{\star}\subset\widehat{M}_n)\geq1-s_n\{b\exp(-cn^a)\}
\end{equation}
where $a, b, c$ are positive values and $s_n$ is the size of true model. The rate $a$ controls how high dimensionality the screening procedure can handle. It will be illustrated in detail in later sections. For a larger $a$, the probability of containing true model converges to 1 faster. The effect of $c$ can be negligible for a larger $a$ and a sufficient large $n$. The constant $b$ is not crucial in the asymptotic sense, but it is important for finite sample situations. With the same value of $a$, a larger value of $b$ indicates that important variables can be miss-selected with higher probability. For some existing results, $b$ can even grow as $n$ increases. Therefore, the constant $b$ and the convergence rate $a$ can be viewed as a measure of effectiveness of a screening method. Correspondingly, a sure screening procedure will be considered more effective with a larger $a$ and a smaller $b$ in the convergence probability of (\ref{j5}). Although the existing screening procedures have been proved to possess the ``sure screening property", the established convergence of containing the true model can be slow subject to various specific model settings and conditions.


Our first goal in the current paper is to develop a general and unified framework for sure screening methods from a loss function perspective. Consider a response variable $Y$ which distribution depends on parameter $\theta$. Suppose $\theta$ is a function of $p-$dimensional covariate vector $\bX=(X_1,\cdots,X_p)^T$. We are interested in selecting the covariates $X_j$'s which are associated with the response variable $Y$ through a nonparametric function $\theta=f(X_1,\cdots,X_p)$. For notational convenience, we will write as $\theta(\bX)$ to denote its dependence on the covariates $\bX$. In later presentation, we sometimes simply write it as $\theta$ for the true function, and the readers should be reminded that the $\theta$ is a function on $\bX$. This setting includes a variety of commonly used regression models:

Example 1 (Gaussian Regression): Assume that $Y|\bX=\bx$ is from $N(\theta(\bx), \sigma^2)$ for some constant $\sigma>0$. 

Example 2 (Logistic Regression): Assume that $Y|\bX=\bx$ is from Bernoulli distribution and  $\ln P(Y=1|\bX=\bx)-\ln P(Y=0|\bX=\bx)=\theta(\bx)$. 

Example 3 (Poisson Regression): Assume that $Y|\bX=\bx$ is from Poisson distribution and $\ln E(Y|\bX=\bx)=\theta(\bx)$.

Example 4 (Quantile Regression): Let $Q_{\alpha}(Y|\bX=\bx)$ be the $\alpha$th quantile of the distribution for $Y|\bX=\bx$, then assume $Q_{\alpha}(Y|\bX=\bx)=\theta(\bx)$. 

The above Examples 1-3 fall within the general framework of mean regression: 
\begin{equation}\label{eq1}
E(Y|\bX=\bx)=h(\theta(\bx))=g^{-1}(f(x_1,\cdots,x_p)),
\end{equation}
where $h$ is some known function, $f$ is a nonparametric function and $g$ is called the link function. When $g$ is the canonical link, that is, $g=(h)^{-1}$, we have $\theta(\bX)=f(X_1,\cdots,X_p)$. However, Example 4 is different from the mean regression. 

The above regression models are equivalent to considering a loss function $l(\omega,Y)$ for measuring the divergence between a generic variable $\omega$ and the response variable $Y$ where $\omega$ is a function of $\bX$, and assuming that the true model of $\theta$ will minimize $E[l(\omega, Y)|\bX=\bx]$ with respect to $\omega$. For instance, in the above Examples 1-3, we can choose $l(\omega, Y)$ as the negative of the log-likelihood of $Y|\bX=\bx$; In the above Example 4, we can choose $l(\omega, Y)=(Y-\omega)[\alpha-\bI(Y-\omega<0)]$, where $\bI$ is an indicator function. Therefore, we will select the important covariates $X_j$'s associated with $Y$ based on such a loss function. 

In the current paper, we newly propose a definition of loss function called conditional strictly convex loss, which contains, but is not limited to, negative log-likelihood loss for one-parameter exponential families, exponential loss for binary classification and quantile regression loss for robust estimation. Our sure screening property is established within such a wide class of loss functions. Therefore, several existing screening methods automatically fall within our framework, including Fan, Feng \& Song (2011) for nonparametric additive models and He, Wang \& Hong (2013) for quantile regression, although their proposed screening procedures can be different from ours. In addition, many more screening methods are suggested by our framework, for example, generalized additive models, binary classification by exponential loss and so on. 

Our second goal of the current paper is to develop screening methods under conditional strictly convex loss with better convergence probability of containing the true model. We treat the marginal regression as fitting the response variable with componentwise covariates via the loss function. We impose an additive model structure for the unknown nonparametric function approximated by B-spline basis. Interestingly, if we consider the goodness of fit statistics as the marginal utility to rank the importance of each covariate to the joint model, we can achieve a much better convergence probability of containing the true model compared with other related literature. Detailed comparison between our results with other related literature will be presented in Section 3. Furthermore, our selected model size can be controlled at the level of sample size $n$ rather than the dimensionality $p_n$ with high probability. 

The major contribution of the current paper is to simultaneously tackle the issues of ``universality" and ``effectiveness". For the ``universality", we establish the sure screening property within a unified framework through the introduction of a new class of loss functions: conditional strictly convex loss; For the ``effectiveness", within this framework, we show that the goodness of fit nonparametric screening methods can achieve a better convergence probability of containing the true model compared with related literature. 

Theoretical pursuit of ``universality" and ``effectiveness" for screening procedures in the current paper has shed new light on the choice of sure screening methods and greatly benefited the applications of screening methods in practice. For example, the superior performance of our proposed method compared with other existing screening procedures will be further demonstrated by extensive simulation studies and some real scientific data example. Our method is called \textbf{G}oodness \textbf{of} \textbf{fi}t \textbf{n}onparametric \textbf{s}creening (Goffins). To stabilize the computation performance, we also provide an iterative screening procedure and an improved variant to handle the situations where covariates are possibly correlated. 

The rest of this paper will be organized as follows: section 2 introduces the conditional strictly convex loss, the B-spline approximation and the goodness of fit nonparametric screening; section 3 establishes the exponential bound, the sure screening properties and the control of model selection size; section 4 proposes an iterative screening procedure and an improved variant; section 5 provides simulation studies and real data analysis. All the technical proofs  and some numerical results are relegated to the supplementary article [Han (2018)].  

\section{Nonparametric Screening under Convex Loss}
\subsection{Conditional Strictly Convex Loss}
Let $l(x,y): \bbR\times\bbR\rightarrow\bbR$ be a function and assume the partial derivative $\partial l(x,y)/\partial x$ exists almost everywhere for $x$ throughout the paper. We consider $l(\omega, Y)$ as a loss function to measure the divergence between a generic variable $\omega$ and the response variable $Y$. We assume the convexity of $l(\omega,Y)$ in the $\omega$ position, that is, $l(t_1\omega_1+t_2\omega_2,Y)\geq t_1l(\omega_1,Y)+t_2l(\omega_2,Y)$ for any real values $t_1+t_2=1$ and $t_1, t_2>0$. Here, $\omega$ is a function of covariates $\bX$, and can be written as $\omega(\bX)$ to denote its dependence on $\bX$. For notational convenience, we sometimes simply write it as $\omega$. Suppose the distribution of $Y$ depends on some parameter $\theta$ where $\theta$ is a nonparametric function of the covariates $\bX$. Assume the true model of $\theta$ minimizes $E[l(\omega,Y)|\bX]$ with respect to $\omega$. 

In the current paper, we will newly propose a definition of loss function called Conditional Strictly Convex Loss. Our sure screening method will be established within such a wide class of convex loss functions.  

\begin{definition}
If $\partial E[l(\omega, Y)|\bX]/\partial \omega$ is continuously differentiable in $\omega$ and $\partial^2 E[l(\omega, Y)|\bX]/\partial \omega^2>0$, then $l(\omega,Y)$ is called a conditional strictly convex loss function. 
\end{definition}

The conditional strictly convex loss includes, but is not limited to, the following three major types of loss functions:

\vspace{0.03in}
\noindent {\bf Type 1: Negative Log-likelihood Loss for Exponential Families} \\
Suppose that the random variable $Y$ is from a one-parameter exponential family with density function 
\begin{equation}\label{h2}
f_{Y|X}(y;\theta)=\exp\big(y\theta-b(\theta)+c(y)\big)
\end{equation}
for some known functions $b()$ and $c()$ where $b''()$ exists. Consider the negative log-likelihood loss:
\begin{equation}\label{h1}
l(\omega,Y)=-[\omega Y-b(\omega)+c(Y)]. 
\end{equation}
Minimization of $E[l(\omega, Y)|\bX]$ with respect to $\omega$ and letting $\theta$ be the minimizer leads to $E[Y|\bX]=b'(\theta)$, which naturally belongs to the mean regression (\ref{eq1}). This is the setting of generalized additive model in Stone (1986). Note that the second derivative of $l(\omega,Y)$ with respect to $\omega$ is $b''(\omega)$, and $b''(\theta)$ is the variance of $Y$ from the exponential families. 

The loss function (\ref{h1}) can be better understood by some popular regression models:  

Example 1 (Gaussian Regression): $b(\theta)=\theta^2/2$, $c(y)=-y^2/2$ and $l(\omega, Y)=(Y-\omega)^2/2$. 

Example 2 (Logistic Regression): $b(\theta)=\ln(1+\exp(\theta))$, $c(y)=0$ and $l(\omega, Y)=-\omega Y+\ln(1+\exp(\omega))$. 

Example 3 (Poisson Regression): $b(\theta)=\exp(\theta)$, $c(y)=-\ln(y!)$ and $l(\omega, Y)=-Y\omega+\exp(\omega)+\ln(Y!)$. \\

\vspace{0.03in}
\noindent{\bf Type 2: Exponential Loss for Classification}\\
In classification problems, suppose $Y\in\{-1,1\}$ and $P(Y=1|\bX=\bx)=p(\bx)$. The goal is to construct a classifier $\theta(\bx)$. When new covariates $\bX$ are available, predict the corresponding class type $Y$ as 1 if $\theta(\bX)>c$ and as -1 if $\theta(\bX)<c$ where $c$ is some threshold. The exponential loss is defined as: 
\begin{equation}\label{h4}
l(\omega, Y)=\exp(-Y\omega), 
\end{equation}
which has been considered as a smooth approximation to the misclassification loss (Freund \& Schapire 1997). Minimization of $E[l(\omega, Y)|\bX]$ with respect to $\omega$ and letting $\theta$ be the minimizer leads to 
\begin{equation*}
\ln \frac{P(Y=1|\bX)}{P(Y=-1|\bX)}=2\theta.
\end{equation*}

\vspace{0.03in}
\noindent {\bf Type 3: Quantile Regression Loss}\\
For many practical problems, the distribution information of response variable $Y$ is usually not available or complicated. Instead of imposing a full distribution, quantile regression framework assumes that the $\alpha$th quantile of $Y$ given $\bX$, $Q_{\alpha}(Y|\bX)$, is some function of $\bX$, thus the distribution assumption can be substantially relaxed (Koenker 2005). Correspondingly, consider the loss function 
\begin{equation}\label{h5}
l(\omega,Y)=(Y-\omega)\{\alpha-\bI(Y-\omega<0)\} 
\end{equation}
for $0<\alpha<1$ where $\bI$ is an indicator function. When $\alpha=1/2$, this is proportional to the least absolute deviation loss $|Y-\omega|$, which is popularly used for robust regression. The loss function $l(\omega,Y)$ is not differentiable in $\omega$. This is a key difference from the aforementioned loss functions. Minimization of $E[l(\omega,Y)|\bX]$ with respect to $\omega$ yields $Q_{\alpha}(Y|\bX)=\theta$ where $\theta$ is the minimizer. 

The following Proposition \ref{t1} shows that with mild conditions, Types 1-3 belong to the conditional strictly convex loss. 
\begin{proposition}\label{t1}
For Type 1, if $b''$ is strictly positive and is a continuous function, then (\ref{h1}) belongs to the conditional strictly convex loss; For Type 2, (\ref{h4}) belongs to the conditional strictly convex loss; For Type 3, if the conditional distribution of $Y|\bX$ has a continuous density function $f_{Y|X}$ and $f_{Y|X}>0$ on any bounded domain, then (\ref{h5}) belongs to the conditional strictly convex loss. 
\end{proposition}
Without any further investigation, one might simply group Types 1 and 2 in Proposition \ref{t1} as one class since the corresponding loss functions are second differentiable in $\omega$. However, we will show in section 3.4 that even for Types 1 and 2, the loss functions possess some fundamental differences in the underlying structures, which raises challenges for proving the model selection size control in section 3.4. 

The name of conditional strictly convex loss is borrowed from ``strictly convex function". However, there are some major differences between the two concepts. If $l(x,y)$ is a strictly convex function in $x$ and $l'(x,y)$ is continuously differentiable in $x$, then $l$ is also a conditional strictly convex loss, but a conditional strictly convex loss might not be a strictly convex function, see Type 3 quantile regression loss as such a counterexample.

A class of convex loss, Bregman divergence, can also be considered here. For a given convex function $q()$ with derivative $q'()$, the Bregman divergence (Bregman 1967) is defined as
\begin{equation}\label{h6}
l(\omega,Y)=q(\omega)-q(Y)+(Y-\omega)q'(Y). 
\end{equation}
Note that $l(\omega, Y)$ is not generally a symmetric function in $\omega$ and $Y$. Suppose $q'()$ is continuously differentiable and $q''()>0$, it is easy to show that such Bregman divergence belongs to the conditional strictly convex loss. It is impossible for us to list all the possibilities here, thus we will not go any further in this direction. It is worth mentioning that the quantile regression loss (\ref{h5}) does not belong to Bregman divergence. More detailed discussions about Bregman divergence are referred to Zhang, Jiang \& Shang (2009).


\subsection{Goodness of Fit Nonparametric Screening}
To capture the nonparametric structure of $\theta(\bX)$, an powerful model for dimensionality reduction is the additive model:
\begin{equation}\label{g12}
\theta(\bX)=m_1(X_1)+\cdots+m_p(X_p)+\mu,
\end{equation}
where $m_j()$ are the square integrable functions and $\mu$ is an unknown constant. For identifiability, we assume $E[m_j(X_j)]=0$ for $j=1,\cdots,p$. Let $M_{\star}=\{j:E[m_j(X_j)]^2>0\}$ be the true sparse model with non-sparsity size $s_n=|M_{\star}|$. Suppose we have observed data $\{(\bX_i,Y_i)\}$ for $i=1,\cdots,n$, which are independent copies of $\{(\bX, Y)\}$. The dimensionality $p$ is ultrahigh and satisfies $\log p=O(n^a)$ for some $a>0$. Based on the sample data, we aim to select a subset of covariates which contains $M_{\star}$ with moderate size. We allow $p$ to grow with $n$, and denote the dimensionality as $p_n$. 

In this paper, we refer to marginal regression as fitting models with componentwise covariates through the loss function $l(\omega, Y)$. We define the population version of the minimizer of the componentwise regression as
\begin{equation}\label{j2}
f_j^M(X_j)\equiv\argmin_{f_j\in L_2(P)}E[l(f_j(X_j),Y)]
\end{equation}
where $P$ denotes the joint distribution of $(\bX, Y)$ and $L_2(P)$ is the class of square integrable functions under measure $P$. We use B-spline basis to approximate the marginal nonparametric regression function. Let $S_n$ be the space of polynomial splines of degree $l\geq1$. Stone (1986) has shown that under some smoothness conditions, the nonparametric functions can be well approximated by functions in $S_n$. Correspondingly, we define
\begin{equation}\label{j9}
f_{nj}^M(X_j)\equiv\argmin_{f_j\in S_n}E[l(f_j(X_j),Y)]. 
\end{equation}
We also define the marginal minimum divergence estimator as
\begin{equation}\label{j10}
\widehat{f}_{nj}^M(X_j)\equiv\argmin_{f_j\in S_n}\bbP_nl(f_j(X_j),Y),
\end{equation}
where $\bbP_ng(\bX,Y)=n^{-1}\sum_{i=1}^ng(\bX_i,Y_i)$ is the empirical expectation for generic function $g()$. Let $\{\Psi_{k}\}_{k=1}^{d_n}$ denote a normalized B-spline basis with $\|\Psi_{k}\|_{\infty}\leq1$, where $\|\cdot\|_{\infty}$ is the sup norm. For any $f_{nj}\in S_n$, we have
\begin{equation}
f_{nj}(x)=\sum_{k=1}^{d_n}\Psi_{k}(x)\beta_{jk}, \ \ 1\leq j\leq p
\end{equation}
for some coefficients $\{\beta_{jk}\}_{k=1}^{d_n}$. The construction of the B spline basis can be found in the well known books, e.g. de Boor (1978). Let $\bPsi_j\equiv\bPsi_j(X_j)=(\Psi_1(X_j),\cdots,\Psi_{d_n}(X_j))^T$, therefore, we can express
\begin{equation}\label{j1}
f_{nj}^M(X_j)=\bPsi_j^T\bbeta_{j}^M, \quad\quad\quad \widehat{f}_{nj}^M(X_j)=\bPsi_j^T\widehat{\bbeta}_{j}^M
\end{equation}
where $\bbeta_{j}^M$ and $\widehat{\bbeta}_{j}^M$ are the $d_n$ dimensional coefficient vector for the minimizers of (\ref{j9}) and (\ref{j10}).

We will consider a sure screening procedure based on goodness of fit statistics. Formally, let 
\begin{equation*}
G_{n,j}=\bbP_n\{l(\widehat{\beta}_0^M,Y)-l(\bPsi_j^T\widehat{\bbeta}_{j}^M,Y)\}, \quad\quad j=1,\cdots,p_n 
\end{equation*}
where $\widehat{\beta}_0^M\equiv\argmin_{\beta_0\in\bbR}\bbP_n l(\beta_0,Y)$. Correspondingly, for the population level,
\begin{equation*}
G_j^{\star}=E\{l(\beta_0^M,Y)-l(\bPsi_j^T\bbeta_{j}^M,Y)\}, \quad\quad j=1,\cdots,p_n
\end{equation*}
where $\beta_0^M\equiv\argmin_{\beta_0\in\bbR}El(\beta_0,Y)$. The goodness of fit statistics compares the marginal regression model with the null model (no variables included in the model). Intuitively, if the marginal contribution of an individual variable is significant to the response variable, the goodness of fit measure should be relatively large. We select model by $\widehat{M}_{\nu_n}=\{1\leq j\leq p_n: G_{n,j}\geq \nu_n\}$ for a predetermined threshold $\nu_n$.  Our screening method is called \textbf{G}oodness \textbf{of} \textbf{fi}t \textbf{n}onparametric \textbf{s}creening (Goffins). We intentionally use the letter ``G" in $G_{n,j}$ and $G_j^{\star}$ to denote the goodness of fit statistics. 

When $l$ is the squared error loss, since the term $\bbP_n l(\widehat{\beta}_0^M,Y)$ in $G_{n,j}$ is not affected by the index $j$, Goffins is equivalent to screening based on the sum of squared residuals, that is, select the model by $\{1\leq j\leq p_n: \bbP_n(Y-\bPsi_j^T\widehat{\bbeta}_{j1}^M)^2\leq \mu_n\}$ for some threshold $\mu_n>0$. Note that, $\bbP_n(Y-\bPsi_j^T\widehat{\bbeta}_{j}^M)^2$ can be further expressed as $\bbP_nY^2-\bbP_n(\bPsi_j^T\widehat{\bbeta}_{j}^M)^2$. Therefore, Goffins under the squared error loss is equivalent to selecting the model by $\{1\leq j\leq p_n: \bbP_n(\bPsi_j^T\widehat{\bbeta}_{j}^M)^2\geq \gamma_n\}$ for some threshold $\gamma_n>0$. More generally, when $l$ is the negative log-likelihood loss for exponential families, Goffins is equivalent to screening based on the likelihood ratio statistics. For parametric model based likelihood ratio screening, see Fan \& Song (2010).   

\section{Sure Screening Properties}
\subsection{Preliminaries}
In this paper, we will show that our goodness of fit nonparametric screening (Goffins) has the sure screening property, and the number of the selected variables has moderate size. Let $[a,b]$ be the support of covariates $X_j$. The following conditions are needed:
\begin{itemize}
\item[(A)] The nonparametric marginal functions $\{f_j^M\}_{j=1}^p$ belong to a class of functions $\mathfrak{F}$ whose $r$th derivative $f^{(r)}$ exists and is Lipschitz of order $\alpha$:
    \begin{equation}
    \mathfrak{F}=\{f(\cdot): |f^{(r)}(s)-f^{(r)}(t)|\leq K|s-t|^{\alpha}, \ \text{for} \ s,t\in [a,b]\}
    \end{equation}
    for some positive constant $K$, where $r$ is a non-negative integer and $\alpha\in(0,1]$ such that $d=r+\alpha>0.5$.
\item[(B)] The marginal density functions $g_j$ of $X_j$ satisfies $0<K_1\leq g_j(X_j)\leq K_2<\infty$ on $[a,b]$ for $1\leq j\leq p$ for some constants $K_1$ and $K_2$.
\item[(C)] The unknown nonparametric function $\theta(\bX)$ satisfies that $\sup_{X\in\mathbb{R}^{p_n}}|\theta(\bX)|<M$ from some positive constant $M$. 
\end{itemize}
Conditions A, B \& C are standard regularity assumptions for nonparametric regression in Stone (1986), Fan, Feng \& Song (2011), He, Wang \& Hong (2013), etc. 

The following Lemma \ref{t2} shows that the approximation error of marginal regression $f_{nj}^M$ in (\ref{j9}) to marginal nonparametric projection $f_j^M$ in (\ref{j2}) is negligible.  
\begin{lemma}\label{t2}
If $l$ is a conditional strictly convex loss, under Conditions A-C, assume that $f_j^M$ is uniformly bounded for $j=1,\cdots,p$, then  there exists a positive constant $C_1$ such that $E(f_j^M-f_{nj}^M)^2\leq C_1d_n^{-2d}$, where $d$ is defined in Condition A.  
\end{lemma}

To show that for $j\in M_{\star}$, $G_j^{\star}$ has a non-vanishing signal, we also need the following conditions:
\begin{itemize}
\item[(D)] $\min_{j\in M_{\star}}E[f_j^M(X_j)-Ef_j^M(X_j)]^2\geq c_1d_nn^{-2\kappa}$, for some $0<\kappa<d/(2d+1)$ and $c_1>0$.
\item[(E)]  $d_n^{-2}\leq c_1(1-\xi)^2n^{-2\kappa}/4C_1$ for some $\xi\in(0,1)$.
\end{itemize}
Condition D requires that the marginal nonparametric projections are at a certain strength level separate from the noise. Therefore, we can select the significant covariates based on a threshold. Similar conditions also appear in related literature on nonparametric screening, e.g. Fan, Feng \& Song (2011) and He, Wang \& Hong (2013). See detailed discussion in section 2 of supplementary article [Han (2018)]. 
\begin{lemma}\label{t5}
Under conditions in Lemma \ref{t2}, in addition, Condition D and E are also satisfied, then  $\min_{j\in M_{\star}}G_j^{\star}\geq \frac{b^{\star}}{2}c_1\xi d_nn^{-2\kappa}$ for some positive constant $b^{\star}$. 
\end{lemma}

As we will show in later sections, the sure screening property depends on the characteristics of a generalized definition for partial derivative of loss function $l(x,y): \mathbb{R}\times\mathbb{R}\rightarrow\mathbb{R}$ with respect to $x$. More specifically, let $\widetilde{l}(x,y):\bbR\times\bbR\rightarrow\bbR$ be a Riemann integrable function with respect to $x$ such that for any $x_1>x_2$ and any $y$, $l(x_1,y)-l(x_2,y)=\int_{x_2}^{x_1}\widetilde{l}(s,y)ds$. Since $l(x,y)$ is differentiable in $x$ almost everywhere, such $\widetilde{l}$ exists and is unique almost everywhere in $x$. For notational convenience, we simply use $l'(x,y)$ to denote one such $\widetilde{l}(x,y)$. When $l(x,y)$ is differentiable in $x$, $l'(x,y)$ is uniquely determined. When we consider the quantile regression loss $l(x,y)=(y-x)\{\alpha-\bI(y-x<0)\}$, if $x>y$, then $l(x,y)=(y-x)(\alpha-1)$; if $x<y$, then $l(x,y)=(y-x)\alpha$. Except at $x=y$, $\partial l(x,y)/\partial x=\bI(y-x<0)-\alpha$. Hence, for quantile regression loss, for any $x_1>x_2$ and any $y$, we have $l(x_1,y)-l(x_2,y)=\int_{x_2}^{x_1}[\bI(y-s<0)-\alpha]ds$. Therefore, we will use $l'(x,y)=\bI(y-x<0)-\alpha$ for the quantile regression loss throughout the paper. The above argument motivates the following Definition 2:  

\begin{definition}
The notation $l'(\omega,Y)$ is defined as follows: for Type 1 and 2, $l'(\omega,Y)=\frac{\partial l(\omega,Y)}{\partial\omega}$; for Type 3, $l'(\omega,Y)=\bI(Y-\omega<0)-\alpha$. 
\end{definition}
To simplify the discussion, for the loss function $l(x,y)$ which is not differentiable in $x$ but is differentiable in $x$ almost everywhere, we only focus on the quantile regression loss here. However, similar argument also applies to other loss functions beyond quantile regression loss. 

To characterize $l'(\omega,Y)$ for the exponential tail bound in section 3.2, we also need the following definition for subgaussian random variables. 
\begin{definition}
A random variable $X$ with mean $\mu=EX$ is called $\sigma-$subgaussian if there is a positive number $\sigma$ such that 
\begin{equation*}
E\exp\big(\lambda(X-\mu)\big)\leq\exp(\frac{\lambda^2\sigma^2}{2}), \quad\quad\quad \forall\lambda\in\mathbb{R}. 
\end{equation*}
\end{definition}
Note that if $X\sim N(\mu,\sigma^2)$, then $X$ is $\sigma-$subgaussian. If random variable $X$ is bounded such that $a\leq X\leq b$, then $X$ is subgaussian with $\sigma=(b-a)/2$. See Buldygin \& Kozachenko (2000) for more details. 


\subsection{Exponential Bound for Marginal Minimum Divergence Estimator}
Since both $\widehat{f}_{nj}^M$ and $f_{nj}^M$ can be expressed in terms of B-spline basis functions, it is crucial to establish an exponential bound for the tail probability of $\|\widehat{\bbeta}_{j}^M-\bbeta_{j}^M\|$. The sharpness of this exponential bound directly affect the convergence probability of the screening method. 

The following Theorem \ref{t3} provides an exponential bound for the tail probability of marginal minimum divergence estimator for the B-spline coefficients. It will serve as the cornerstone for our later derivations of the other theorems. The following conditions are required for Theorem \ref{t3}. 
\begin{itemize}
\item[(F)] $d_n=o(n^{1/3})$ and $d_n=O(n^{2\kappa})$. 
\item[(G)] $E[l'(\omega,Y)|\bX]$ is bounded for any bounded $\omega$.  
\end{itemize}
\begin{proposition}
For Types 1-3, under the conditions in Proposition \ref{t1}, condition G is satisfied. 
\end{proposition}

The tail probability of $\|\widehat{\bbeta}_{j}^M-\bbeta_{j}^M\|$ depends on the properties of $l'(\omega,Y)$. More specifically, we will consider the following set of conditions:
\begin{itemize}
\item[(H1)] $l'(\omega,Y)$ is bounded for any bounded $\omega$;
\item[(H2)] $l'(\omega,Y)$ conditional on $\bX$ is a $\sigma-$subgaussian random variable where $\sigma$ does not depend on $\bX$; 
\item[(H3)] For any bounded $\omega$, $E[\exp\big(\lambda l'(\omega,Y)\big)|\bX]<\infty$ for all $|\lambda|\leq c_0$ with some constant $c_0>0$, 
\end{itemize}
The notation $l'(\omega,Y)$ in Condition G and H1-H3 is based on Definition 2. By Definition 3, if Condition H1 is satisfied, then Condition H2 is also satisfied; if Condition H2 is satisfied, then Condition H3 is also satisfied. In the following Theorem \ref{t3}, we will show that with stronger assumption a better tail probability can be correspondingly achieved. 

To better understand the wide applicability of Conditions H1-H3, let us consider some examples from Types 1-3 which satisfy these conditions. Some popular regression models can be summarized in the following Proposition \ref{t4}. 

\begin{proposition}\label{t4}
Types 2 and 3 satisfy Condition H1. For Type 1, if $Y|\bX$ follows Bernoulli distribution, then (\ref{h1}) satisfies Condition H1; if $Y|\bX$ follows Normal distribution, then (\ref{h1}) satisfies Condition H2; if $Y|\bX$ follows Poisson distribution, then (\ref{h1}) satisfies Condition H3. 
\end{proposition}

Furthermore, if $Y|\bX$ follows some other distributions in the exponential family, under some regularity conditions, it is possible that the corresponding loss function (\ref{h1}) also satisfy Condition H3. For example, if $Y|\bX\sim \text{Laplace}(\mu(\bX),b)$ with a known parameter $b$, then Condition H3 is satisfied. If $Y|\bX\sim \text{Exponential}(\lambda(\bX))$, and if there exists a positive constant $c$ such that $\lambda(\bX)\geq c$, then Condition H3 is satisfied. Similar arguments for verifying Condition H3 also apply to Chi-square distribution, negative binomial distribution,  inverse-Gaussian distribution with a known shape parameter, Gamma distribution with a known scale parameter. To save space, we will not discuss in detail for these examples. 

\begin{theorem}\label{t3}
For a convex loss $l(\omega,Y)$, if it is also a conditional strictly convex loss, for any constant $c_3>0$, under Conditions C, F, G, 
there exists positive constants $c_4$ and $c_5$ such that for sufficiently large $n$, \\
if Condition H1 is satisfied and $d_n=o(n^{1-2\kappa})$, then 
\begin{equation}\label{h16}
P(\|\widehat{\bbeta}_{j}^M-\bbeta_{j}^M\|^2\geq c_3d_nn^{-2\kappa})\leq\exp(-c_4n^{1-2\kappa}d_n^{-1});
\end{equation}
if Condition H2 is satisfied and $d_n=o(n^{(1-2\kappa)/2})$, then 
\begin{equation}\label{s6}
P(\|\widehat{\bbeta}_{j}^M-\bbeta_{j}^M\|^2\geq c_3d_nn^{-2\kappa})\leq\exp(-c_4n^{1-2\kappa}d_n^{-1})+\exp(-c_5n^{1-2\kappa}d_n^{-2});
\end{equation}
if Condition H3 is satisfied and $d_n=o(n^{(1-2\kappa)/3})$, then 
\begin{equation}\label{s7}
P(\|\widehat{\bbeta}_{j}^M-\bbeta_{j}^M\|^2\geq c_3d_nn^{-2\kappa})\leq\exp(-c_4n^{1-2\kappa}d_n^{-1})+2\exp(-c_5n^{1/2-\kappa}d_n^{-3/2}).
\end{equation}
\end{theorem}
In (\ref{h16}), when $d_n=o(n^{1-2\kappa})$, $n^{1-2\kappa}d_n^{-1}$ in the tail probability diverges to infinity as $n$ increases, which implies that the tail probability converges to zero. Similar arguments also apply to (\ref{s6}) and (\ref{s7}). It is worth mentioning that Theorem \ref{t3} is proved based on a unified argument with some modifications according to each situation of Conditions H1-H3. The proof is different from the related literature and can be of independent research interest. 

\subsection{Sure Screening}
Based on Theorem \ref{t3} for estimation of B-spline coefficients, we are now ready to establish the sure screening property for our Goffins method. Different properties of loss functions can lead to different convergence probabilities of containing the true model. 
\begin{theorem}\label{t6}
Under the conditions in Theorem \ref{t3} and Lemma \ref{t5}, \\
(i) for Types 1 and 2, there exists a positive constant $\zeta$, then by taking $\nu_n=\nu d_nn^{-2\kappa}$ with $0<\nu\leq\zeta$, there exists positive constants $c_4$, $c_5$ and $c_6$ such that
if Condition H1 is satisfied and $d_n=o(n^{1-2\kappa})$, then
\begin{equation}\label{s8}
P(M_{\star}\subset\widehat{M}_{\nu_n})\geq1-s_n\Big[\exp(-c_4n^{1-2\kappa}d_n^{-1})+6\exp(-c_5n^{1-2\kappa})\Big];
\end{equation}
if Condition H2 is satisfied and $d_n=o(n^{(1-2\kappa)/2})$, then
\begin{eqnarray}\label{s9}
P(M_{\star}\subset\widehat{M}_{\nu_n})&\geq&1-s_n\Big[\exp(-c_4n^{1-2\kappa}d_n^{-1})\\
&&+\exp(-c_5n^{1-2\kappa}d_n^{-2})+6\exp(-c_6n^{1-2\kappa})\Big];\nonumber
\end{eqnarray}
if Condition H3 is satisfied and $d_n=o(n^{(1-2\kappa)/3})$, then
\begin{eqnarray}\label{s10}
P(M_{\star}\subset\widehat{M}_{\nu_n})&\geq&1-s_n\Big[\exp(-c_4n^{1-2\kappa}d_n^{-1})\\
   &&+2\exp(-c_5n^{1/2-\kappa}d_n^{-3/2})+6\exp(-c_6n^{1-2\kappa})\Big];\nonumber
\end{eqnarray}
(ii) for Type 3, if $d_n=o(n^{1-2\kappa})$, take $\nu_n=\nu d_nn^{-2\kappa}$ with $\nu\leq b^{\star}c_1\xi/4$ where $b^{\star}$ is defined in Lemma \ref{t5}, there exists positive constants $c_4$ and $c_5$ such that 
\begin{equation}\label{s11}
P(M_{\star}\subset\widehat{M}_{\nu_n})\geq1-s_n\Big[\exp(-c_4n^{1-2\kappa}d_n^{-1})+12\exp(-c_5n^{1-2\kappa})\Big]. 
\end{equation}
\end{theorem}
Theorem \ref{t6} shows that our Goffins method corresponding to conditional strictly convex loss possesses the sure screening property. It follows from Theorem \ref{t6} that in (\ref{s8}) and (\ref{s11}) we can handle the NP-dimensionality: $\log p_n=o(n^{1-2\kappa}d_n^{-1})$. Under this condition, $P(M_{\star}\subset\widehat{M}_{\nu_n})\rightarrow1$ to achieve the sure screening property. For (\ref{s9}), the NP-dimensionality will be changed to $\log p_n=o(n^{1-2\kappa}d_n^{-2})$ and for (\ref{s10}) we can handle $\log p_n=o(n^{1/2-\kappa}d_n^{-3/2})$. 

The proof of Theorem \ref{t6} in the supplementary article [Han (2018)] is not limited to Types 1-3. For example, let Class A be the loss functions such that $l''(\omega,Y)\equiv\partial^2l(\omega,Y)/\partial\omega^2$ exists, $l''(\omega, Y)$ is continuous in $\omega$, $l''(\omega, Y)>0$ and $l''(\omega,Y)$ is bounded when $\omega$ is bounded, then the results corresponding to Types 1-2 in Theorem \ref{t6} are also valid for the loss functions in Class A. It is not difficult to verify that Types 1-2 are only special examples in Class A. When $l(\omega, Y)$ is not differentiable in $\omega$, the discussion is more complicated. Let $\widetilde{l}(x,y):\bbR\times\bbR\rightarrow\bbR$ be a Riemann integrable function with respect to $x$ such that for any $x_1>x_2$ and any $y$, $l(x_1,y)-l(x_2,y)=\int_{x_2}^{x_1}\widetilde{l}(s,y)ds$. Since we assume that the loss function $l(x,y)$ is differentiable in $x$ almost everywhere, such $\widetilde{l}$ exists and is unique almost everywhere in $x$. Let Class B be the loss functions such that there exists a corresponding $\widetilde{l}$ where $\widetilde{l}(\omega,Y)$ is bounded and $\widetilde{l}(\omega,Y)$ is non-decreasing in $\omega$, then the results corresponding to Type 3 is also valid for the loss functions in Class B. Our definition of $l'(\omega,Y)$ for quantile regression loss clearly satisfies the conditions in Class B. 

\subsection{Controlling Selection Size}
The sure screening methods will not be informative unless the model selection size can be controlled at a reasonable level. The following Theorem \ref{t7} shows that our Goffins method can control the size of the selected variables at the level of the sample size $n$ rather than the dimension $p_n$. However, for controlling the model size, our definition of conditional strictly convex loss is not sufficient for the discussion. We need a finer class of loss functions which possesses certain structures. This is the motivation of our following Definition 4. 
\begin{definition}
If a convex loss $l(\omega,Y)$ satisfies that 
\begin{equation*}
\partial E[l(\omega,Y)|\bX]/\partial\omega=G(\omega)-H(\bX)K(\omega)
\end{equation*}
for some functions $G( )$, $H( )$ and $K( )$, then $l(\omega,Y)$ is called conditional derivative separable loss. 
\end{definition}
\begin{proposition}
Types 1-3 are conditional derivative separable loss functions. 
\end{proposition}
\noindent For Type 1, $\partial E[l(\omega, Y)|\bX]/\partial\omega=b'(\omega)-\mu(\bX)$ where $\mu(\bX)=E(Y|\bX)$. Therefore, $G(\omega)=b'(\omega)$, $H(\bX)=\mu(\bX)$ and $K(\omega)=1$. 

\vspace{0.03in}
\noindent For Type 2, $\partial E[l(\omega, Y)|\bX]/\partial\omega=\exp(\omega)-[\exp(\omega)+\exp(-\omega)]p(\bX)$ where $p(\bX)=E(Y|\bX)$. Therefore, $G(\omega)=\exp(\omega)$, $H(\bX)=p(\bX)$ and $K(\omega)=\exp(\omega)+\exp(-\omega)$. 

\vspace{0.03in}
\noindent For Type 3, $\partial E[l(\omega, Y)|\bX]/\partial\omega=F_{Y|X}(\omega)-\alpha$ where $F_{Y|X}$ is the conditional cumulative distribution function of $Y|\bX$. Therefore, $G(\omega)=F_{Y|X}(\omega)$, $H(\bX)=\alpha$ and $K(\omega)=1$. 

\vspace{0.03in}
Detailed discussions reveal different structures of loss functions. For example, Types 1 and 3 both have a constant $K(\omega)=1$ while Type 2 does not have such property. Furthermore, Types 1 and 2 are second differentiable in $\omega$ but Type 3 is not differentiable in $\omega$. Fortunately, we can propose a unified proof to bound $\sum_{j=1}^{p_n}E(f_{nj}^M-Ef_{nj}^M)^2$ when the loss function $l(\omega, Y)$ is conditional derivative separable loss and conditional strictly convex loss, which will serve as a major step for controlling the model selection size.
\begin{theorem}\label{t7}
Let $\bPsi=(\bPsi_1,\cdots,\bPsi_{p_n})^T$, $\bbeta^{\star}$ be the coefficient vector of basis functions for the joint regression model of $\theta(\bX)$ on $\bX$, $\beta_0^{\star}$ be the intercept term in the joint regression model and $\bSigma=E\bPsi\bPsi^T$. If $l$ is a conditional derivative separable loss and conditional strictly convex loss, under conditions in Theorem \ref{t6}, in addition, $E(\bPsi^T\bbeta^{\star})^2=O(1)$ and $K(\beta_0^{\star})\neq0$, then we have \\
(i) $\sum_{j=1}^{p_n}E(f_{nj}^M-Ef_{nj}^M)^2=O(d_n\lambda_{\max}(\bSigma))$; \\
(ii) with $\nu_n$ described in Theorem \ref{t6}, there exist constants $c_4$, $c_5$, $c_6$ such that\\
\noindent\text{Types 1 and 2:} if Condition H1 is satisfied and $d_n=o(n^{1-2\kappa})$, then 
\begin{equation*}
P\Big(|\widehat{M}_{\nu_n}|\leq O\big(n^{2\kappa}\lambda_{\max}(\bSigma)\big)\Big)\geq1-p_n\Big[\exp(-c_4n^{1-2\kappa}d_n^{-1})+6\exp(-c_5n^{1-2\kappa})\Big];
\end{equation*}
if Condition H2 is satisfied and $d_n=o(n^{(1-2\kappa)/2})$, then
\begin{eqnarray*}
P\Big(|\widehat{M}_{\nu_n}|\leq O\big(n^{2\kappa}\lambda_{\max}(\bSigma)\big)\Big)&\geq&1-p_n\Big[\exp(-c_4n^{1-2\kappa}d_n^{-1})\\
&&+\exp(-c_5n^{1-2\kappa}d_n^{-2})+6\exp(-c_6n^{1-2\kappa})\Big];
\end{eqnarray*}
if Condition H3 is satisfied and $d_n=o(n^{(1-2\kappa)/3})$, then
\begin{eqnarray*}
P\Big(|\widehat{M}_{\nu_n}|\leq O\big(n^{2\kappa}\lambda_{\max}(\bSigma)\big)\Big)&\geq&1-p_n\Big[\exp(-c_4n^{1-2\kappa}d_n^{-1})\\
&&+2\exp(-c_5n^{1/2-\kappa}d_n^{-3/2})+6\exp(-c_6n^{1-2\kappa})\Big];
\end{eqnarray*}
\text{Type 3:} if $d_n=o(n^{1-2\kappa})$, then 
\begin{equation*}
P\Big(|\widehat{M}_{\nu_n}|\leq O\big(n^{2\kappa}\lambda_{\max}(\bSigma)\big)\Big)\geq1-p_n\Big[\exp(-c_4n^{1-2\kappa}d_n^{-1})+12\exp(-c_5n^{1-2\kappa})\Big].
\end{equation*}
\end{theorem}

The tail probabilities directly follow from Theorem \ref{t6} and have been explained in Section 3.3. The selected model size depends on the matrix $\bSigma$ which involves the dependence structure of the covariates. As discussed in Fan, Feng \& Song (2011), it can be assumed that $\lambda_{\max}(\bSigma)=n^{\tau}$ for some $\tau>0$. Correspondingly, the model size will be controlled at a reasonable rate of $n$. With iterative Goffins method described in the next section 4, we will show in the simulation studies that the number of false positives can be very small while the true important variables are all selected even when the covariates are correlated. 


\subsection{Connection and Comparison with Related Literature}
When $l$ is the squared error loss, Fan, Feng \& Song (2011)'s screening for nonparametric additive models is based on $n^{-1}\sum_{i=1}^n(\widehat{f}_{nj}^M(X_{i,j}))^2$. They presented a similar result to Theorem \ref{t6} here but with a different convergence probability as
\begin{equation*}
P(M_{\star}\subset\widehat{M}_{\nu_n})\geq 1-s_nd_n\Big\{(8+2d_n)\exp(-c_4^{*}n^{1-4\kappa}d_n^{-3})+6d_n\exp(-c_5^{*}nd_n^{-3})\Big\}, 
\end{equation*}
and they can handle the NP dimensionality $\log p_n=o(n^{1-4\kappa}d_n^{-3})$. Compared with their result, our convergence probability in (\ref{s9}) for gaussian regression not only improve on the convergence rate, but also improve significantly on those coefficient terms. We can achieve NP-dimensionality $\log p_n=o(n^{1-2\kappa}d_n^{-2})$. It should be noted that Fan, Feng \& Song (2011) established the result under a weaker assumption than our Condition H2. More specifically, they assume $Y=\theta(\bX)+\epsilon$, $E(\epsilon|\bX)=0$ and for any $B_1>0$, $E[\exp(B_1|\epsilon|)|\bX]\leq B_2$ for some constant $B_2$. On the other hand, this condition is stronger than our Condition H3. If we use (\ref{s10}) for the comparison here in favor of Fan, Feng \& Song (2011)'s result, then we can handle the NP-dimensionality $\log p_n=o(n^{1/2-\kappa}d_n^{-3/2})$. When $n^{1/3-2\kappa}=O(d_n)$, we can handle a higher dimensionality. Otherwise, their result is better. 

When $l$ is the quantile regression loss, He, Wang \& Hong (2013)'s screening is based on $n^{-1}\sum_{i=1}^n(\widehat{f}_{nj}^M(X_{i,j}))^2$. They have shown that for positive constants $c_6^*$, $c_7^*$ and $c_8^*$, 
\begin{equation*}
P\Big(\|\widehat{\bbeta}_{j}^M-\bbeta_{j}^M\|^2\geq c_6^*d_nn^{-4\kappa}\Big)\leq 2\exp(-c_7^*n^{1-8\kappa})+\exp(-c_8^*n^{1-4\kappa}d_n^{-2}). 
\end{equation*}
Correspondingly, they presented a convergence probability as
\begin{equation*}
P(M_{\star}\subset\widehat{M}_{\nu_n})\geq1- s_n\Big\{11\exp(-c_9^{*}n^{1-8\kappa})+12d_n^2\exp(-c_{10}^{*}n^{1-4\kappa}d_n^{-3})\Big\}. 
\end{equation*}
Note that He, Wang \& Hong (2013) considers a signal strength in Condition D as $c_1n^{-2\kappa}$, and the parameter $\tau$ in Theorem 3.3 of He, Wang \& Hong (2013) is equivalent to $2\kappa$ in our paper here. If we reset the minimum signal strength in Condition D the same as that of He, Wang \& Hong (2013), our result for Theorem \ref{t3} will be modified as 
\begin{equation*}
P(\|\widehat{\bbeta}_{j}^M-\bbeta_{j}^M\|^2\geq c_3d_nn^{-4\kappa})\leq\exp(-c_4n^{1-4\kappa}d_n^{-1})
\end{equation*}
under the conditions of He, Wang \& Hong (2013). Correspondingly, our result for Theorem \ref{t6} will be modified as
\begin{equation*}
P(M_{\star}\subset\widehat{M}_{\nu_n})\geq1-s_n[\exp(-c_4n^{1-4\kappa}d_n^{-1})+12\exp(-c_5n^{1-4\kappa})]. 
\end{equation*}
Therefore, He, Wang \& Hong (2013)'s convergence probability also indicates larger possibility of not selecting important covariates for a non asymptotic setting. 

When $l$ is the negative log-likelihood loss for one-parameter exponential families, Fan \& Song (2010) constructed sure screening for generalized linear models. If the B-spline approximation is treated as a type of group variable selection, then our Theorem \ref{t3} has some connections with Fan \& Song (2010)'s result. Compared with Fan \& Song (2010), the tail probability in our Theorem \ref{t3} does not have the extra term $nP(\Omega_n^c)$ in their paper where $n$ is the sample size and $\Omega_n$ is the region such that the loss function satisfies some Lipschitz condition. In Fan \& Song (2010), their exponential bound also involves a Lipschitz constant. When the response variable is not bounded (e.g., most of the exponential families), this Lipschitz constant diverges to infinity, which results in a slower convergence rate for the tail bound, in contrast with our result. For example, when considering the squared error loss, Fan \& Song (2010) Theorem 4 will have 
\begin{equation*}
P(\|\widehat{\bbeta}_{j}^M-\bbeta_{j}^M\|\geq c_3n^{-\kappa})\leq \exp(-c_4n^{(1-2\kappa)/3})+nm_1\exp(-c_4n^{(1-2\kappa)/3}).
\end{equation*}
for bounded covariates. For our Theorem \ref{t3} (under Condition H2), let $d_n=2$ for a fair comparison, then we have 
\begin{equation*}
P(\|\widehat{\bbeta}_{j}^M-\bbeta_{j}^M\|\geq c_3n^{-\kappa})\leq 2\exp(-c_4n^{1-2\kappa}). 
\end{equation*}
It is clear that we have a much better result here. When considering the Poisson regression loss, the corresponding convergence rate in the tail probability bound can be much slower than $(1-2\kappa)/3$. For our Theorem \ref{t3} (under Condition H3), let $d_n=2$ for a fair comparison, we have 
\begin{equation*}
P(\|\widehat{\bbeta}_{j}^M-\bbeta_{j}^M\|\geq c_3n^{-\kappa})\leq \exp(-c_4n^{1-2\kappa})+2\exp(-c_5n^{1/2-\kappa}). 
\end{equation*}
It is still a better result than Fan \& Song (2010). 

\section{Iterative Goffins Method and Improved Variant}
In practice, unimportant variables can be correlated with the important variables, therefore such variables can have significant marginal effects even though they are not significant in the joint true model. To improve the performance of our screening method, we consider an iterative version of Goffins. Given the data $\{(\bX_i, Y_i)\}$, $i=1,\cdots,n$, we choose the same truncation term $d_n=O(n^{1/5})$. In Theorem \ref{t6}, the threshold $\nu_n$ is chosen at the level $d_nn^{-2\kappa}$. In practice, the parameter $\kappa$ is unknown, but we can determine a data-driven threshold. To achieve this, we extend the random permutation idea of Fan, Feng \& Song (2011) and Zhao \& Li (2012). Let $\bX$ be the matrix with the $i$th row as $\bX_i$. The algorithm works as follows:
\begin{itemize}
\item[Step 1.] For every $j\in\{1,\cdots,p\}$, compute 
\begin{equation*}
\widehat{f}_{nj}=\argmin_{ f_{nj}\in S_n}\bbP_nl(f_{nj}(X_j), Y)\quad\quad1\leq j\leq p.
\end{equation*}
Randomly permute the rows of $\bX$ and we have $\widetilde{\bX}=(\widetilde{X}_1,\cdots,\widetilde{X}_p)$. Let $\omega_{(q)}$ be the $q$th quantile of $\{ G_{n,j}^*, j=1,\cdots,p\}$, where $\widehat{f}_{nj}^*=\argmin_{f_{nj}\in S_n}\bbP_nl(f_{nj}(\widetilde{X}_j), Y)$. Then our method selects the following variables: $\mathcal{A}_1=\{j: G_{n,j}\geq \omega_{(q)}\}$. In our numerical studies, we choose $q=1$, the maximum value of the empirical norm of the permuted estimates. 
\item[Step 2.] Apply penalized regression on the set $\mathcal{A}_1$ to select a subset $\mathcal{M}_1$. Specifically, when $l$ is the negative log-likelihood loss, apply the penalized generalized additive model regression (e.g. penGAM in Meier, van de Geer \& B\"uhlmann 2009). 
\item[Step 3.] For every $j\in\mathcal{M}_1^c=\{1,\cdots,p\}/\mathcal{M}_1$, minimize $\bbP_nl(f_0+\sum_{i\in\mathcal{M}_1}f_{ni}(X_i)+f_{nj}(X_j), Y)$ with respect to $f_0\in\mathbb{R}$, $f_{ni}\in S_n$ for all $i\in\mathcal{M}_1$ and $f_{nj}\in S_n$. For identifiability, we apply the B spline basis without the intercept for $j\in \mathcal{M}_1^c$ and for $i\in\mathcal{M}_1$. Apply the screening procedure with adaptive threshold determined by the new random permutation. Choose a set of indices $\mathcal{A}_2$. Then penalized regression is applied on the set $\mathcal{M}_1\bigcup\mathcal{A}_2$ to select a subset $\mathcal{M}_2$. 
\item[Step 4.] Iterate the process until $|\mathcal{M}_l|\geq s_0$ or $\mathcal{M}_l=\mathcal{M}_{l-1}$. 
\end{itemize}

This iterative version of Goffins will be denoted as ``I-Goffins" in our simulation studies. To further stabilize the performance, we can apply a ``cap" to control the number of selected variables in each iteration. For example, in our simulation studies, we restrict to select 1 variable at each step. Since the chance of selecting unimportant variables in each step has been reduced, the probability of selecting important variables in the subsequent steps has been improved. This is the idea behind the greedy INIS method proposed by Fan, Feng \& Song (2011) for additive modeling. To be consistent with Fan, Feng \& Song (2011), we name this improved variant of our method as greedy iterative goodness-of-fit nonparametric screening (GI-Goffins). 

\section{Simulation Studies}
Similar to Fan, Feng \& Song (2011), we set $n=400$ but we consider $p=1000, 2000, 5000$ for all examples to investigate the impact of high dimensionality on screening methods. Following Fan, Feng \& Song (2011), we consider the number of spline basis functions as $d_n=\lceil n^{1/5}\rceil+2=6$. Note that in this paper we consider the full B spline basis, and Fan, Feng \& Song (2011) considered the B spline basis without the intercept. The goodness of fit screening methods under the two sets of basis are equivalent. Eight simulation examples will be constructed according to the four major types of regressions: Gaussian regression, Logistic regression, Poisson regression and quantile regression.  Define 
\begin{eqnarray*}
&&f_1(x)=x,\quad f_2(x)=(2x-1)^2,\quad f_3(x)=\sin(2\pi x)/(2-\sin(2\pi x)),\\
&&f_4(x)=0.1\sin(2\pi x)+0.2\cos(2\pi x)+0.3\sin(2\pi x)^2\\
&&\quad\quad\quad+0.4\cos(2\pi x)^3+0.5\sin(2\pi x)^3,\\
&&f_5(x)=\exp(x-0.5),\quad f_6(x)=0.1\sin(2\pi x)^2+0.4\cos(2\pi x)^3\\
&&f_7(x)=\sin(x-1),\quad f_8(x)=(x-1.5)^2,\quad f_9(x)=2\cos(x)/(2-\sin(x)). 
\end{eqnarray*}
\begin{itemize}
\item Model 1 (Linear Regression): $Y|\bX=5f_1(X_1)+3f_2(X_2)+4f_3(X_3)+6f_4(X_4)+\sqrt{1.74}\epsilon$. Each $X_i\sim Uniform(0,1)$ i.i.d. and $\epsilon\sim N(0,1)$. 
\item Model 2 (Linear Regression): The model is the same as Model 1 but the covariates $\bX=(X_1,\cdots, X_p)^T$ are simulated according to the random effects model $X_j=(W_j+tU)/(1+t), j=1,\cdots,p$ where $W_1,\cdots, W_p$ and $U$ are i.i.d. $Uniform(0,1)$ and $t=0.4$. 
\item Model 3 (Logistic Regression): $\ln\big(P(Y=1|\bX)/P(Y=0|\bX)\big)=2f_1(X_1)+3f_7(X_2)+2f_8(X_3)+3.5f_9(X_4)$. Each $X_i\sim Uniform(-2.5,2.5)$ i.i.d.
\item Model 4 (Logistic Regression): The model is the same as Model 3 but the covariates $\bX=(X_1,\cdots, X_p)^T$ are simulated according to the random effects model $X_j=(W_j+tU)/(1+t), j=1,\cdots,p$ where $W_1,\cdots, W_p$ are i.i.d. from $Uniform(-2.5,2.5)$, independent of $U\sim Uniform(0,1)$ and $t=0.4$. 
\item Model 5 (Poisson Regression): $Y|\bX\sim Poisson\big(\exp\{f_1(X_1)+f_3(X_2)+f_5(X_3)+f_6(X_4)\}\big)$. Each $X_i\sim Uniform(0,1)$ i.i.d.
\item Model 6 (Poisson Regression): The model is the same as Model 5 and the covariates $\bX=(X_1,\cdots, X_p)^T$ are simulated according to the same structure as Model 2. 
\item Model 7 (Heteroscedastic Regression): $Y|\bX=5f_1(X_1)+3f_2(X_2)+4f_3(X_3)+4f_5(X_4)+0.5\exp(f_6(X_{20})+f_7(X_{21})+f_8(X_{22}))\epsilon$, where $\bX\sim N_p(0,\bSigma)$ independent of $\epsilon\sim N(0,1)$ and the $(i,j)$th element of covariance matrix $\bSigma$ is $0.8^{|i-j|}$. 
\item Model 8 (Heteroscedastic Regression): The model is the same as Model 7 except that the random error $\epsilon\sim Laplace(0,2)$.  
\end{itemize}
Models 1 and 2 have been similarly considered in Meier, van de Geer \& B\"uhlmann (2009) and Fan, Feng \& Song (2011), while Models 3-8 are newly proposed in the current paper. The covariates are independent in Models 1, 3, 5 but correlated in Models 2, 4, 6, 7, 8. Note that Model 8 is different from Model 7 because the Laplace(0,2) distribution for random error will emphasize more on the covariates $X_{20}$, $X_{21}$, $X_{22}$, thus making the heteroscedatic regression model more challenging.

\begin{table}
\caption{The median (IQR) and $5\%$, $25\%$, $75\%$, $95\%$ quantiles of minimum model size.  }
\begin{center}
\begin{tabular}{c|ccccccc}
\hline\hline
        Model   &$p$   & Methods   & Median (IQR)  & $5\%$  &$25\%$   &$75\%$   &  $95\%$    \\
\hline                       
       Model  3   &1000     &Goffins  &  5(7.25)     &4               &   4              &  11.25      &67.05   \\
                         &              & Kfilter   &14(28)         &  4             &  6              &  34         &  115.05   \\
                          &             &QaSIS         &    NA           &   NA      &   NA    &  NA   & NA   \\
                          &             &SIS         &   393.5(537)     &  38.9       &  158     & 695     & 941.05     \\
                        &               &SIRS        &  393(537)         & 38.95             &   158              &   695     & 941.05  \\
                        &               &DC        &   16(25)                  &   4          &  9                 &  34                & 112.05  \\
                        &               &EL         &  451(524.25)       & 55.75             &  211                &  725.25        & 948.05   \\  \cline{2-8}                
                         & 2000    &Goffins  &  5(6)                &  4             &  4               &   10                &  64.05 \\
                         &             & Kfilter   &  13(28.25)                & 4              &  7              &     35.25              & 119.65    \\
                         &              &QaSIS         &   NA               &    NA         &    NA              &     NA             &  NA  \\
   Size=4                       &             &SIS         & 376.5(521)                 & 29.85              &   168.50              &  689.5                & 954.05    \\
                          &             &SIRS        &  376(521.5)                  &29.85              &  168               &  689.5                 &  954.05 \\
                         &              &DC        &  16(24.25)                   & 4            &   8                &  32.25                & 105.20  \\
                         &              &EL         & 441(508.25)                    &  38.95            &  215.50                &  723.75               &  961.05  \\  \cline{2-8}                  
                         & 5000     &Goffins  & 9(36.25)                 &  4             &  4               &   40.25                & 301.20   \\
                         &             & Kfilter   & 59(151.5)                 & 5              &  20              &   171.5                &  597   \\
                         &              &QaSIS         &  NA                & NA              &  NA                &   NA               &NA    \\
                          &             &SIS         & 2164(2819.5)         &165.1               & 887.5                & 3707.0                 & 4814.4   \\
                          &             &SIRS        & 2163.5(2817.5)         & 164.20             &  888.75               &3706.25                   &4814.40   \\
                         &              &DC        &67(130.5)                     & 8            &   27.75                & 158.25                 &504.15   \\
                         &              &EL         &  2465.5(2705.75)                   & 229.85             & 1136.25                &  3842.00               & 4846.05   \\                                                              
\hline   
       Model  4    &1000       &Goffins  & 5(6)                 & 4              &  4               &  10                 &55.1   \\
                         &           & Kfilter   & 16(36)                 &   4            &  6              & 42                  & 163.1    \\
                          &             &QaSIS         &   NA               &  NA             &  NA                &  NA                & NA   \\
                          &             &SIS         &   418.5(566)               &  20.9             &  148.75               &  714.75         & 920.3   \\
                        &               &SIRS        &    418.5(566.25)          & 19.95         & 148.50                &   714.75        &920.3   \\
                        &               &DC        &  17(33)                   & 4            & 8                  &  41                & 137.05  \\
                        &               &EL         &   495(543.75)        & 42             & 211.75                 & 755.50                 &  933.15  \\  \cline{2-8}
                        &2000              &Goffins  &  6(17)                & 4               &    4             &  21                 &131.1   \\
                         &          & Kfilter   &   32.5(72.25)               & 4              &10                &    82.25               & 243.35    \\
                         &              &QaSIS         &   NA               &    NA           &     NA             &  NA                &NA    \\
   Size=4                       &             &SIS         &  858(1079.25)          & 36.85              & 312.25                &1391.50        & 1892   \\
                          &             &SIRS        &  858(1077)              & 36.9             &  314.5               &  1391.5           & 1892  \\
                         &              &DC        &   34(70)                  &  5           &   14                &   84               &  240.15 \\
                         &              &EL         &   1008(1037.25)           & 83.85             &  446.75         & 1484         & 1911.05    \\    \cline{2-8} 
                         &5000        &Goffins  & 9(32)                 &  4             & 4                &    36               &332.4   \\
                         &             & Kfilter   &  59.5(174)                &   5            & 16               &190                   &  770.75   \\
                         &              &QaSIS         &    NA              &   NA            &    NA              &   NA               &NA    \\
                          &             &SIS         &  1857.5(2778.5)                & 84.95              & 749.50                &   3528.00               &  4675.10  \\
                          &             &SIRS        & 1855.5(2781.5)                   & 84.95             & 746.50                 &  3528.00                 & 4675.10  \\
                         &              &DC        &  70(160.5)                   & 8             & 26                  &   186.5               & 661.0   \\
                         &              &EL         &  2261.5(2701)                   & 196.8              & 1054.5                 & 3755.5                & 4726.5   \\     
\hline\hline
\end{tabular}
\end{center}
\end{table}

\begin{table}
\caption{The median (IQR) and $5\%$, $25\%$, $75\%$, $95\%$ quantiles of minimum model size.  }
\begin{center}
\begin{tabular}{c|ccccccc}
\hline\hline
        Model   &$p$   & Methods   & Median (IQR)  & $5\%$  &$25\%$   &$75\%$   &  $95\%$    \\
\hline                       
       Model  5    &1000              &Goffins  &   4(0)               &    4           &   4              &       4            &  12.1 \\
                         &            & Kfilter   &  28(63)                &  4             &    11            &    74               &  265.05   \\
                          &             &QaSIS         &  NA                &   NA            & NA                 &    NA            & NA   \\
                          &             &SIS         &   447(525.75)               &    22.95           & 174.25                &     700             &   944.15 \\
                        &               &SIRS        &  491(501)                  & 46             & 242.5                &    743.5               &  951.05  \\
                        &               &DC        &  11(18)                   &   4          &  6                 &   24               & 71.05   \\
                        &               &EL         &   461.5(510)                  & 32.9             &     194.75             &  704.75          & 949   \\  \cline{2-8}
                         &2000              &Goffins  &  4(1)                & 4              &   4              &    5               & 17  \\
                         &            & Kfilter   &   60(144)               & 5              &   17             & 161                  & 504.1    \\
                         &              &QaSIS         &    NA              &  NA             &      NA            &     NA             &NA    \\
   Size=4                       &             &SIS         &  980(1154.25)          & 48.95              &   371.75              &   1526               & 1924.05    \\
                          &             &SIRS        & 1031(1028)              & 112.75             & 516                &  1544                 & 1904.30   \\
                         &              &DC        & 21(36)                    & 4             &     9              &  45                &  132 \\
                         &              &EL         &  1013(1166.75)           & 54.9             & 393.75                 & 1560.5                & 1923.05   \\     \cline{2-8}
                          &5000              &Goffins  & 4(2)                 &4               &     4            &     6            &  42 \\
                         &      & Kfilter   &  147.5(377)                & 7              & 41.75               &   418.75         & 1340.00    \\
                         &              &QaSIS         &  NA                &  NA             &  NA                &  NA                &  NA  \\
                          &             &SIS         & 2237.5(2662.75)      &163.80               & 978.25                & 3641.00           & 4748.45   \\
                          &             &SIRS        &  2619.5(2518.75)                  & 266.80             & 1468.75               &  3987.50          & 4723.20  \\
                         &              &DC        &  44(98.25)                   & 6            &   17                &  115.25         & 360.70  \\
                         &              &EL         &  2320(2701)                   & 183.85             & 988.00                 &  3689.00         & 4773.15   \\                                                              
\hline   
       Model  6   & 1000             &Goffins  &   4(1.25)               &4               &4                 &      5.25             & 19.05   \\
                         &               & Kfilter   &  35.5(95)                &   4.95            & 12               &   107       & 345.05    \\
                          &             &QaSIS         &  NA                & NA              &       NA           &   NA               &  NA  \\
                          &             &SIS         &   400.5(571.75)        & 28               &  150.25               &  722        & 940.05   \\
                        &               &SIRS        &  510.5(514)             & 49.95             &  244.5               &  758.5        &952.15   \\
                        &               &DC        &  17(36)                   &  4           &  7                 &  43                &  129 \\
                        &               &EL         &   425(567)                  &  30.95            &  163                & 730                & 940.05    \\  \cline{2-8}
                         &2000              &Goffins  & 4(2)                 & 4              &   4              &  6                 & 31   \\
                         &             & Kfilter   & 66.5(159)                 & 5              & 22               &   181          & 552.7    \\
                         &              &QaSIS         &    NA              &  NA             &   NA               &     NA             & NA   \\
   Size=4                       &             &SIS         &  930.5(1067.25)                &  41.95             &  356                & 1423.25                  &  1905.10  \\
                          &             &SIRS        &  1027.5(1051.75)            & 94.80             & 488.50                &  1540.25             & 1911.05  \\
                         &              &DC        &  23(47)                   &5             &   10             &  57                & 240.05  \\
                         &              &EL         & 972.5(1059.5)       & 50.80             & 383.25                 & 1442.75                & 1912   \\  \cline{2-8}                        
                         &5000      &Goffins  &  5(7)                &  4             &   4              &      11             &   74.4 \\
                         &            & Kfilter   & 182(483.25)                 & 6              &   52             &    535.25               &  1614.10    \\
                         &              &QaSIS         &  NA                &     NA          &    NA              &       NA           & NA   \\
                          &             &SIS         &  2072.5(2690.75)           & 198              &825                 &     3515.75             &4709.70    \\
                          &             &SIRS        & 2442(2528.25)              & 231.70             & 1324                &  3852.25                 &  4778.45 \\
                         &              &DC        &  74(173.5)                   & 6             &    23.75               &  197.25                & 593.15  \\
                         &              &EL         &   2130(2627.5)                  &   214.6           & 926.5                 &  3554.0                & 4710.3   \\     
\hline\hline
\end{tabular}
\end{center}
\end{table}

\noindent \textbf{Minimum Model Size} Following Fan \& Song (2010), Fan, Feng \& Song (2011) as well as later literature in sure screening field, we use the minimum model size required to contain the true model $M_{\star}$ as a measure of the effectiveness of a screening method. The simulation round is 500 for all the examples. We compare our Goffins method with six other successful screening methods in the existing literature, including some recent model-free screening methods. More specifically, we consider fused Kolmogorov filter (Kfilter) by Mai \& Zou (2015), quantile adaptive screening (QaSIS) by He, Wang \& Hong (2013), SIS for generalized linear model by Fan \& Song (2010), sure independent ranking and screening (SIRS) by Zhu, Li, Li \& Zhu (2012), distance correlation learning (DC) by Li, Zhong \& Zhu (2012) and empirical likelihood screening (EL) by Chang, Tang \& Wu (2013). Note that when considering squared error loss, our Goffins is equivalent to NIS by Fan, Feng \& Song (2011). Therefore, we will treat NIS as a special example of Goffins, and will not present NIS as a separate method for comparison here. However, when considering quantile regression loss, our Goffins method is different from QaSIS in He, Wang \& Hong (2013), because Goffins is based on goodness of fit statistics while QaSIS is based on squared norm of fitted nonparametric function. We choose quantile $75\%$ whenever a quantile regression loss is considered. The implementations of QaSIS and SIRS are based on http://users.stat.umn.edu/$\sim$wangx346/research/example1b.txt. The implementations of Kfilter, DC and EL are based on the R codes from the authors of related literature. The implementation of SIS is based on the R package ``SIS". 

In Tables 1-2 along with the Tables S1-S2 in the supplementary article [Han (2018)], we present the median, the interquartile range (IQR) and different quantiles of minimum model size. Following existing literature, if the median is closer to the true model size and the IQR is smaller, the corresponding screening method is considered as more effective. Overall, our Goffins method performs best among the seven screening methods. For Models 1, 3, 4, 5, 6, our medians of minimum model size are close to the true model size 4 and IQRs are the smallest. Our medians and IQRs will not increase significantly when the dimensionality $p$ increases from 1000 to 5000. For comparison, other methods tend to select a much larger model to contain the true model, and the performance can deteriorate dramatically when $p$ increases. Furthermore, $5\%$, $25\%$, $75\%$ and $95\%$ quantiles of our minimum model size are significantly smaller than the other methods. Models 7-8 are very challenging heteroscedastic regression models, but our method still performs better than the other methods, including Kfilter and QaSIS. Table S2 in the supplementary article [Han (2018)] also suggests that even when we consider quantile regression loss, Goffins is different from QaSIS. Model 2 turns out to be a difficult example for all the methods. However, our simulation in tables S3-S5 of supplementary article [Han (2018)] will show that an iterative version of Goffins (GI-Goffins or I-Goffins) can substantially reduce the false positives while selecting the true important variables. 

\section{Data Analysis}
Classification between the malignant pleural mesothelioma (MPM) and the lung cancer adenocarcinoma (ADCA) has received increasing attention in both clinical studies and high dimensional statistical research. Gordon, et al (2002) studied the data from 181 tissue samples (31 MPM and 151 ADCA) with 12533 gene expression levels for each sample. Among these 181 sample data, 16 MPM and 16 ADCA have been combined as the training set while the other 149 samples (15 MPM and 134 ADCA) are considered as the testing set. The goal of research is in two-fold as explained in Gordon, et al (2002): 1. Find the minimum number of predictor genes that are most importantly associated with the disease type; 2. Construct a classifier rule which can predict the future patients' disease type based on their gene expression levels with high statistical accuracy. Aspect 1 can substantially reduce the medical cost of obtaining patients' relevant gene data and the cost of potential scientific experiments on such genes. The performance of the classifier is usually evaluated based on the testing data. 

Since the disease type is a categorical data, and the number of genes is extremely high ($p=12533$) compared with the small sample size ($n=32$), we will apply our GI-Goffins method with respect to the logistic regression for the training data. We first standardize the gene expression data for each gene over the training samples such that the sample mean is 0 and the sample standard deviation is 1. Our method selects five genes that are importantly associated with the disease type: ``31575-f-at", ``37716-at", ``39795-at", ``41286-at" and ``41402-at". We construct a generalized additive model (B spline basis without the intercept and the number of spline basis functions as $d_n=\lceil n^{1/5}\rceil+1=3$) based on such five genes and apply the model to the training data. The fitted nonparametric functions corresponding to those five genes have also been plotted in Figure S1 in the supplementary article [Han (2018)]. Then we apply our constructed model to the test data. Among the 149 samples for the testing data, we make 144 correct predictions. For the 5 samples that we misclassified, one MPM sample has been predicted as ADCA while four ADCA samples has been predicted as MPM. ISIS for the generalized linear model has also been considered to select important variables. To be fair, we also apply a generalized additive model based on the selected genes for the training data and further use this fitted model for classification on the test data. However, this method will select fewer and different genes and the performance is much inferior to our method. I-EL is an iterative version of EL and penalized empirical likelihood regression described in Chang, Tang \& Wu (2013). Its performance is even worse than ISIS.


This lung cancer data has also been analyzed by various statistical methods in the past literature. It is impossible and unnecessary for us to list all the relevant results here, and we only compare our method with some representative methods which have been shown superior performance. In Table 3 , we will compare our GI-Goffins method with linear discriminant methods such as ROAD in Fan, Feng \& Tong (2012) and FAIR in Fan \& Fan (2008). Our GI-Goffins is a good balance between the testing error and the number of selected genes compared with other methods. More selected genes will cause substantial cost in future diagnosis and experiments. Therefore, GI-Goffins is the method that we recommend for practice. 

\begin{table}\label{t8}
\caption{Performance of methods on lung cancer data. $p=12533$. }
\begin{center}
\begin{tabular}{cccc}
\hline\hline
        Method   & Training Error & Testing Error  & Number of Selected Genes   \\
\hline                       
        GI-Goffins        & 0/32          &  5/149       &   5           \\
        ISIS              & 0/32          &  18/149     &    2          \\
        I-EL                &  0/32         &  40/149     &   5          \\
\hline        
        ROAD          & 1/32         &  1/149         &  52               \\
        FAIR            &  0/32        &   7/149       &   31          \\                 
\hline\hline
\end{tabular}
\end{center}
\end{table}

\section{Further Discussions}
\subsection{Optimality}
An interesting question is whether the convergence rate in the upper bound of the  tail probability that we established in Theorem \ref{t3} is optimal. More specifically, if we have
\begin{equation*}
b_1\exp(-c_2n^a)\leq P(\|\widehat{\bbeta}_{j}^M-\bbeta_{j}^M\|^2\geq c_1d_nn^{-2\kappa})\leq b_2\exp(-c_3n^a)
\end{equation*}
for some constants $a$, $b_1$, $b_2$, $c_1$, $c_2$ and $c_3$, then we can say that the convergence rate $a$ in the upper bound of the tail probability is optimal, because the convergence rate $a$ can not be improved further. 

When the loss function $l$ is the negative log likelihood loss of one-parameter exponential families, under general regularity conditions, the maximum likelihood estimator has the asymptotic normality (Heyde 1997, Gao, et al. 2008), that is, 
\begin{equation*}
[I_j(\bbeta_{j}^M)]^{1/2}(\widehat{\bbeta}_{j}^M-\bbeta_{j}^M)-N(0,\bI_{d_n})\rightarrow 0\quad\quad\text{in distribution}
\end{equation*}
where $I_j$ is the information matrix of the $j$th covariate. Plugging in the negative log likelihood loss and the B-spline basis functions, we have
\begin{equation*}
n^{1/2}\{E[b''(\bPsi_j^T\bbeta_{j}^M)\bPsi_j\bPsi_j^T]\}^{1/2}(\widehat{\bbeta}_{j}^M-\bbeta_{j}^M)-N(0,\bI_{d_n})\rightarrow0\quad\quad\text{in distribution}
\end{equation*}
Since $\bPsi_j^T\bbeta_{j}^M$ is bounded based on the argument in the supplementary article [Han (2018)] and $b''()$ is a continuous function, $b''(\bPsi_j^T\bbeta_{j}^M)$ is upper bounded by a positive constant. Furthermore, due to Lemma 3 in the supplementary article [Han (2018)], 
\begin{equation*}
(\widehat{\bbeta}_{j}^M-\bbeta_{j}^M)^Tn\{E[b''(\bPsi_j^T\bbeta_{j}^M)\bPsi_j\bPsi_j^T\}(\widehat{\bbeta}_{j}^M-\bbeta_{j}^M)\leq D_2nd_n^{-1}\|\widehat{\bbeta}_{j}^M-\bbeta_{j}^M\|^2. 
\end{equation*}
Therefore, asymptotically, we have 
\begin{eqnarray*}
P(\|\widehat{\bbeta}_j^M-\bbeta_j^M\|^2\geq c_1d_nn^{-2\kappa})&=&P(D_2nd_n^{-1}\|\widehat{\bbeta}_j^M-\bbeta_j^M\|^2\geq D_2c_1n^{1-2\kappa})\\
              &=&P(\chi_{d_n}^2\geq D_2c_1n^{1-2\kappa}). 
\end{eqnarray*}
Thus, we need to find a lower bound for the tail probability of $\chi_{d_n}^2$ distribution. When $d_n=1$, it is well known that for any positive $y$, 
\begin{equation*}
P(\chi_1^2\geq y)\geq 1-\sqrt{1-\exp(-\frac{2y}{\pi})}=\frac{\exp(-\frac{2y}{\pi})}{1+\sqrt{1-\exp(-\frac{2y}{\pi})}}\geq \frac{1}{2}\exp(-\frac{2y}{\pi}). 
\end{equation*}
Let $y=D_2c_1n^{1-2\kappa}$, comparing with our Theorem \ref{t3} under Condition H1 or H2, we have achieved the optimal convergence rate $n^{1-2\kappa}$ asymptotically. When $d_n=2$, for any positive $y$, $P(\chi_2^2\geq y)=\exp(-\frac{y}{2})$. Comparing with our Theorem \ref{t3} under Condition H1 or H2, we have also achieved the optimal convergence rate $n^{1-2\kappa}$ asymptotically. For more general $d_n$, we do not have a sharp lower bound of the tail probability of Chi-square distribution. Therefore, we will not discuss further here. 

\subsection{Adaptive Threshold}
Theorem \ref{t6} is established based on a threshold $\nu_n$ at the level of $d_nn^{-2\kappa}$. In practice, the parameter $\kappa$ is unknown. Therefore, we need an adaptive threshold for the real data. Consider a threshold $\widehat{\nu}_n$ which is constructed based on the sample data, it will be interesting to derive a lower bound for $P(M_{\star}\subseteq\widehat{M}_{\widehat{\nu}_n})$. We have
\begin{equation*}
P(M_{\star}\subseteq\widehat{M}_{\widehat{\nu}_n})\geq 1-\sum_{j\in M_{\star}}P(G_{n,j}<\widehat{\nu}_n). 
\end{equation*}
Note that 
\begin{equation*}
P(G_{n,j}<\widehat{\nu}_n)\leq P(G_{n,j}<\nu_n)+P(\widehat{\nu}_n>\nu_n). 
\end{equation*}
We have derived the upper bound of $P(G_{n,j}<\nu_n)$ is the proof of Theorem \ref{t6}. Therefore, we need to derive an upper bound for the second term here. 

Consider a permutation of the sample covariates $\{\bX_i\}_{i=1}^n$. We can obtain the estimates of marginal regression based on the permuted data:
\begin{equation*}
(\widehat{\bbeta}_{j}^M)^{\pi}=\argmin_{\beta_j\in \bbR^{d_n}} \bbP_nl(\bPsi_j^T(X_j^{\pi})\bbeta_j,Y)
\end{equation*}
where $\pi=(\pi_1,\cdots,\pi_n)$ is a permutation of the index $\{1, 2, \cdots, n\}$. Note that $(G_{n,j})^{\pi}$ is a statistical estimate of 0. We will derive an upper bound for $P(\widehat{\nu}_n>\nu)$ for a special case where the loss function $l()$ is the squared error loss. Note that the least squares estimate follows
\begin{equation*}
\sqrt{n}\widehat{\bbeta}_{j}^{\pi}=(\frac{1}{n}\bPsi(X_j^{\pi})\bPsi(X_j^{\pi})^T)^{-1}(\frac{1}{\sqrt{n}}\sum_{i=1}^n\bPsi^T(X_j^{\pi_i})Y_i). 
\end{equation*}
By Anderson \& Robinson (2001) Theorem \ref{t6}, after some algebra, we have the asymptotic normality: 
\begin{equation*}
\sqrt{n}\widehat{\bbeta}_{j}^{\pi}-N(0, \widetilde{\bSigma})\rightarrow 0\quad\quad\text{in distribution}, 
\end{equation*}
where $\widetilde{\bSigma}=(E\bPsi\bPsi^T)^{-1}[E(\bPsi-E\bPsi)(\bPsi-E\bPsi)^T](E\bPsi\bPsi^T)^{-1}$. 

If we let $\widehat{\nu}_{n,j}=\frac{1}{n}\sum_{i=1}^n(\bPsi^T(X_j^{\pi_i})\widehat{\bbeta}_{j}^{\pi})^2$, different thresholds for the marginal utilities of different covariates. As we have discussed in the section 2.2, this screening based on $\bbP_n(\bPsi_j^T\widehat{\bbeta}_{j}^M)^2$ is equivalent to our goodness of fit screening when the loss function is the squared error loss. 
We can show that asymptotically
\begin{equation*}
P(\widehat{\nu}_{n,j}\geq c_1d_nn^{-2\kappa})\leq P(\chi_{d_n}^2\geq c_2d_nn^{-2\kappa})\leq \exp(-c_3d_nn^{1-2\kappa})
\end{equation*}
for some constants $c_1$, $c_2$ and $c_3$. The second inequality is by Laurent \& Massart (2000) Lemma \ref{t2}. Correspondingly, we have a lower bound for the convergence probability of containing the true model. In practice, the threshold $\widehat{\nu}_n$ will be chosen as the maximum value of $\{\widehat{\nu}_{n,j}\}_{j=1}^{p_n}$ under a number of permutations, and the loss function can be general. We do not have a theoretical result for such more complicated situations. 

\subsection{Choice of Loss Function}
The framework of goodness of fit nonparametric screening includes many screening methods based on the choice of loss functions. An important question is how to choose loss function for practical data. When the response variable $Y$ takes values $\{0,1\}$, we suggest to consider the logistic regression loss: $l(\omega, Y)=-\omega Y+\ln(1+\exp(\omega))$; When $Y$ takes nonnegative integer values, we suggest to consider the Poisson regression loss: $l(\omega, Y)=-Y\omega+\exp(\omega)+\ln(Y!)$; When the distribution of $Y$ is expected to be complicated, the quantile regression loss can be considered; When $Y$ is a continuous variable, we suggest to start with the Gaussian regression loss: $l(\omega, Y)=(Y-\omega)^2/2$. This brief guideline could raise misspecification issue of loss functions. 

\subsection{Iterative Screening Procedure}
The idea of iterative screening and penalization has been proposed since Fan \& Lv (2008), and has achieved numerical success in practice. However, formal theoretical justification is still an open problem in the field. The first step is a marginal screening. To simplify the discussion, assume a fixed threshold $\gamma_n$ is applied the selected variables $\mathrm{A}_1=\{j: G_{n,j}\geq\gamma_n\}$ satisfies $M_{\star}\subseteq\mathrm{A}_1$ with high probability (sure screening property). For the second step, based on the set $\mathrm{A}_1$, we apply some penalized regression and select a subset $\mathrm{M}_1$. Ideally, we want to show sign consistency for $\mathrm{M}_1$ under some regularity conditions. The difficulty is that the set $\mathrm{A}_1$ is random, which is different from a conventional penalization regression. Fortunately, we can borrow the technique in Weng, Feng \& Qiao (2017), which considers a two-step procedure for linear regression model (similar to screening $+$ penalization). For the third step, it is a conditional marginal screening after penalization. Emre, Fan \& Verhasselt (2016) has shown the sure screening property based on the conditional screening for generalized linear model. Therefore, if the sign consistency is achieved in step 2, then under some regularity conditions, sure screening property can be achieved in step 3. By mathematical induction, the iterative procedure can achieve sign consistency. We would like to explore the technical details as our future studies. 

\vspace{0.2in}
\noindent\textbf{Acknowledgements}\\
The author wants to thank the Joint Editor, the Associate Editor and the three anonymous referees for many insightful comments which significantly improve the presentation of the paper. 

The author deeply appreciates Professor Jianqing Fan for his encouragement and constructive comments on this project. The author would like to thank Dr. Cheng Yong Tang and Professor Linda Zhao for helpful discussions on the paper. The author also thanks Dr. Shujie Ma and Dr. Lily Wang for the helpful discussion on B spline basis. 

Special thanks go to the following researchers who kindly share their codes of numerical studies: Dr. Yang Feng for NIS, Dr. Qing Mai for Kfilter, Dr. Lukas Meier for penGAM, Dr. Cheng Yong Tang for EL and Dr. Wei Zhong for DC.   

\begin{supplement}[id=suppA]
  \sname{Supplement A}
  \stitle{Supplement to ``Nonparametric Screening under Conditional Strictly Convex Loss for Ultrahigh Dimensional Sparse Data"}
  \slink[doi]{COMPLETED BY THE TYPESETTER}
  \sdatatype{.pdf}
  \sdescription{Due to the space limit, all the technical proofs as well as some numerical results are relegated to the supplementary article [Han (2018)]. }
\end{supplement}

\end{document}